\documentclass[11pt,fullpage]{article}
\usepackage{graphicx}
\usepackage{latexsym}
\usepackage{color}

\usepackage{amsmath}
\usepackage{graphicx}
\usepackage{amssymb}
\usepackage{amsthm}
\usepackage{amsfonts}

\usepackage{eufrak}

\long\def\ignore#1{}

\renewcommand{\baselinestretch} {1.3}

        \makeatletter 
\def\singlespace{\def\baselinestretch{1}\@normalsize}

\newtheorem{theorem}{Theorem}
\newtheorem{lemma}{Lemma}

\newtheorem{corollary}{Corollary}
\newtheorem{proposition}{Proposition}

\newcommand{\be}{\begin{equation}}
\newcommand{\ee}{\end{equation}}
\newcommand{\beqn}{\begin{eqnarray}}
\newcommand{\eeqn}{\end{eqnarray}}
\newcommand{\bt}{\begin{theorem}}
\newcommand{\et}{\end{theorem}}
\newcommand{\bl}{\begin{lemma}}
\newcommand{\bp}{\begin{proposition}}
\newcommand{\ep}{\end{proposition}}
\newcommand{\bc}{\begin{corollary}}
\newcommand{\ec}{\end{corollary}}

\newcommand{\bfm}[1]{\mbox{\boldmath $#1$}}
\newcommand{\bbeta}{\bfm{\beta}}
\newcommand{\bepsilon}{\bfm{\epsilon}}
\newcommand{\by}{{\bf y}}
\newcommand{\bd}{{\bf d}}
\newcommand{\bu}{{\bf u}}
\newcommand{\bg}{{\bf g}}
\newcommand{\bw}{{\bf w}}
\newcommand{\Mp}{\mathfrak{M}_{p_0}}

\begin{document}
\title{\Large{\bf MODEL SELECTION IN REGRESSION UNDER STRUCTURAL CONSTRAINTS}}

\author{
{\bf Felix Abramovich}\\
Department of Statistics and Operations Research\\
Tel Aviv University\\
Tel Aviv 69978, Israel\\
felix@post.tau.ac.il
\\ \\
{\bf Vadim Grinshtein}\\
Department of Mathematics and Computer Science\\
The Open University of Israel\\
Raanana 43537, Israel\\
vadimg@openu.ac.il}

\date{}
\maketitle

\begin{abstract}
The paper considers model selection in regression under the additional
structural constraints on admissible models where the number of potential
predictors might be even larger than the available sample size.
We develop a Bayesian formalism which is used as a natural
tool for generating a wide class of model selection criteria based on
penalized least squares estimation with
various complexity penalties associated with a prior on a model size.
The resulting criteria are
adaptive to structural constraints. We establish the upper bound for the
quadratic risk of the resulting MAP estimator and the corresponding lower bound for
the minimax risk over a set of admissible models of a given size. We then
specify the class of priors (and, therefore, the class of complexity penalties)
where for the ``nearly-orthogonal'' design the MAP estimator is
asymptotically at least nearly-minimax (up to a log-factor) simultaneously
over an entire range of sparse and dense setups. Moreover, when the numbers
of admissible models are ``small'' (e.g., ordered variable selection) or,
on the opposite, for the case of complete variable selection,
the proposed estimator achieves the exact minimax rates.
\end{abstract}



\section{Introduction} \label{sec:intr}
Consider the standard Gaussian linear regression model
\begin{equation}
\by=X\bbeta+\bepsilon, \label{eq:model}
\end{equation}
where ${\bf y} \in \mathbb{R}^n$ is a vector of the observed response variable
$Y$, $X_{n \times p}$ is the design matrix of the $p$ explanatory variables
(predictors)  $X_1,\ldots,X_p$,
$\bbeta \in \mathbb{R}^p$ is a vector of unknown regression coefficients,
$\bepsilon \sim N({\bf 0},\sigma^2I_n)$ and the noise variance $\sigma^2$ is assumed to be known.

A variety of statistical applications of regression
models in different fields nowadays involves a vast number of potential
predictors.
Moreover, $p$ might
be even large relative to the amount of available data $n$ ($p \gg n$ setup)
that raises a severe ``curse of dimensionality'' problem. However, typically
only some of the predictors have a truly relevant impact on the response
$\by$. Model (variable) selection by identifying the ``best'' sparse subset
of these ``significant'' predictors becomes therefore crucial
in the analysis of such large data sets.
For a selected model (a subset of predictors) $M$, the corresponding coefficients
$\bbeta_M$ are then typically estimated by least squares.

The goodness of model selection
depends on the particular goal at hand.
One should distinguish, for example, between estimation of regression
coefficients $\bbeta$, estimation of the mean vector $X\bbeta$, model
identification and predicting future observations.  Different goals may lead to
different optimal model selection procedures especially when the number of
potential predictors $p$ might be much larger than the sample size $n$.
In this paper we consider mainly the estimation of the
mean vector $X\bbeta$ and the goodness of a model
$M$ is measured by the quadratic risk $E||X\hat{\bbeta}_M-X\bbeta||^2$, where
$\hat{\bbeta}_M$ is the least squares estimate of $\bbeta$ for $M$.
The ``best'' model then is the one with the minimal quadratic risk.
Note that the true underlying model in (\ref{eq:model}) is not necessarily
the best in this sense since sometimes it is possible to reduce its
risk by excluding predictors with small (but still nonzero!) coefficients.

Minimum quadratic risk criterion for model selection is evidently impossible to implement since
it involves the unknown true $\bbeta$ but the corresponding ideal minimal
(oracle) risk
can be used as a benchmark for any available model selection procedure.
Typical model selection criteria are
based on minimizing the {\em empirical} quadratic risk
$||{\bf y}-X\hat{\bbeta}_M||^2$, which is the least squares,
penalized by a complexity penalty $Pen(|M|)$ increasing with a model size
$|M|$:
\begin{equation}
\min_M \left\{||{\bf y}-X\hat{\bbeta}_M||^2+Pen(|M|) \right\}  \label{eq:pmle}
\end{equation}

The properties of the resulting penalized least squares estimator
depend obviously on the particular choice of the complexity penalty
$Pen(\cdot)$ in (\ref{eq:pmle}). The most commonly used choice is a
{\em linear} type penalty of the form $Pen(k)=2\sigma^2 \lambda k$
for some $\lambda>0$. The most known examples include $C_p$
(Mallows, 1973) and AIC (Akaike, 1973) for $\lambda=1$, BIC
(Schwarz, 1978) for $\lambda=(\ln n)/2$ and RIC (Foster \& George,
1994) for $\lambda=\ln p$. A series of recent works proposed the
so-called $2k\ln(p/k)$-type {\em nonlinear} complexity penalties of
the form $Pen(k)=2\sigma^2\lambda k (\ln(p/k)+1)(1+o(1))$ with
$\lambda \geq 1$. (see, e.g., Birg\'e \& Massart, 2001, 2007;
Abramovich {\em et al.}, 2006; Bunea, Tsybakov \& Wegkamp, 2007;
Abramovich \& Grinshtein, 2010; Rigollet \& Tsybakov, 2011).

As we have mentioned, in the analysis of large complex data sets it is typically
believed that the underlying (unknown) model is sparse, where a natural measure
of model's sparsity is its size $p_0$.
Abramovich \& Grinshtein (2010) and Rigollet \& Tsybakov (2011)
established the minimax rates for the quadratic risk
of estimating the mean vector $X\bbeta$ in (\ref{eq:model}) over
the set of models of sizes at most $p_0$
(see also analogous results of Raskutti, Wainwright \& Yu, 2012 for the
case when $X$ is of a full rank, i.e. $rank(X)=\min(p,n)$).
They also showed that for $2k\ln(p/k)$-type
penalties, the resulting penalized estimators (\ref{eq:pmle}) are simultaneously
asymptotically optimal
(in the minimax sense) for the entire range of sparse and
 dense models, while linear penalties cannot achieve such a wide optimality
range.

So far, the above minimax properties of various model selection
procedures in (\ref{eq:model}) have been established for
{\em complete variable selection}, where the set of admissible models
contains all $2^p$ subsets of the predictors $X_1,...,X_p$. However,
in a variety of regression setups there exist additional {\em
structural constraints} that restrict the set of admissible models.
In some cases predictors have some natural order and $X_j$ can
enter the model only if $X_{j-1}$ is already there ({\em ordered
variable selection}). For example, in polynomial regression, higher
order polynomials are usually considered only if polynomials of
lower degrees are already in the model. In tree-based models, a
predictor cannot enter a model unless all its ancestors are there
({\em hierarchical model selection}). Another important example of
hierarchical model selection is
models with interactions, where interactions are typically selected only with
the corresponding main effects. For a model with factor predictors,
each factor predictor with $k$ levels is, in fact, associated with a
group of $k-1$ indicator (dummy) predictors. In this case either none or all
of this group can be selected. It is somewhat similar to group sparse
models, where predictors are splitted into pre-defined groups (see, e.g.
Lounici {\em at al.}, 2011).

In this paper we investigate
the minimaxity properties of model
selection criteria in (\ref{eq:pmle}) over classes of sparse and
dense models under additional general structural constraints
extending the existing results for the complete variable selection.
The key point is to adapt the choice of the complexity penalty in (\ref{eq:pmle})
to the specific structural constraints.
In fact, it turns out that it is only the number of admissible models of a given
size that matters.

We extend a Bayesian approach to model selection developed in Abramovich \&
Grinshtein (2010) for the complete variable selection. The proposed
Bayesian formalism is based on imposing a prior on a model size,
where the penalty term in (\ref{eq:pmle}) is then naturally treated
as proportional to its logarithm. From the Bayesian perspective, the
model selection criterion (\ref{eq:pmle}) corresponds thus to the
maximum {\em a posteriori} (MAP) Bayes rule. Depending on a specific
choice of a prior, it implies a variety of penalized least squares
estimators (or, equivalently, model selection criteria) with
different complexity penalties. We show that under mild conditions
on the prior, the resulting MAP estimator falls within a general
class of penalized least squares estimators (\ref{eq:pmle})
considered in Birg\'e \& Massart (2001, 2007) with complexity
penalties $Pen(|M|)$ satisfying certain technical conditions depending
on the number of admissible models of size $|M|$, that allows us to
derive the upper bound for its quadratic risk. On the other hand, we
establish the corresponding lower bound for the minimax quadratic
risk over a set of models of a given size under structural
constraints. We then specify the class of priors (and, therefore,
the class of the corresponding complexity penalties), where for the
``nearly-orthogonal'' design, the resulting MAP estimators
asymptotically are at least nearly-minimax (up to a log-factor)
simultaneously over the wide range of sparse and dense models. In
particular, when the numbers of admissible models are ``small'' (e.g.,
ordered variable selection) and for complete variable selection
they lead to AIC-type and $2k\ln(p/k)$-type estimators
respectively and achieve the exact minimax rate.

The paper is organized as follows. In Section \ref{sec:map} we
develop a Bayesian formalism for model selection in regression under
structural constraints and derive the MAP estimator. Its theoretical
properties are presented in Section \ref{sec:main}. In particular,
we establish the upper bound for its quadratic risk, the
corresponding lower bound for the minimax risk and discuss its
asymptotic minimaxity in various setups. Section \ref{sec:simul}
presents the results of a simulation study. Concluding remarks are
summarized in Section \ref{sec:tam}, while all the proofs are given
in the Appendix.

\section{MAP model selection procedure under structural constraints} \label{sec:map}
We first extend the Bayesian formalism for the model selection in
linear regression developed by Abramovich \& Grinshtein (2010) to
structural constraints. We assume that the latter are known and the
set of all admissible models is therefore fixed.

Consider the linear regression model (\ref{eq:model}), where
the number of possible predictors $p$ might be even larger then the number of
observations $n$. Let $r=rank(X) (\leq \min(p,n))$ and assume that any
$r$ columns of $X$ are linearly independent.
For the ``standard'' linear regression setup, where all $p$
predictors are linearly independent and there are at least $p$ linearly
independent design points, $r=p$.

Any model $M$ is uniquely  defined by the $p \times p$ diagonal indicator matrix
$D_M=diag({\bf d}_M)$, where $d_{Mj}=\mathbb{I}\{X_j \in M\}$ and, therefore,
$|M|=tr(D_M)$. For a given model $M$, we estimate its coefficients by least
square estimator
$\hat{\bbeta}_M=(D_M X'X D_M)^+ D_M X' \bf{y}$,
where ``+'' denotes the generalized inverse matrix.

Let $m(p_0)$ be the number of all admissible models of size $p_0$.
The case $p_0=0$ corresponds to a null model with a single intercept
and, therefore, $m(0)=1$. In fact, we can consider only $p_0 \leq r$
since otherwise, there necessarily exists another vector $\bbeta^*$
with at most $r$ nonzero entries such that $X\bbeta=X\bbeta^*$.
Although for $p_0=r$ there may be several different admissible
models, all of them are evidently undistinguishable for estimating
$X\bbeta$ and can be associated with a single (saturated) model.
Thus, without loss of generality, we can always assume that
$m(r)=1$. Obviously, for any $p_0$, $0 \leq m(p_0) \leq {p \choose
p_0}$, where the two extreme cases $m(p_0)=1$ and $m(p_0)={p \choose
p_0}$ for all $p_0=0,\ldots,r-1$, correspond respectively to the
ordered and complete variable selection.

We start from imposing a prior on the model size $\pi(k)=P(|M|=k),\;k=0,\ldots,r$.
Obviously, $\pi(k)=0$ and $P(M \bigl|\;|M|=k)=0$ iff there are no admissible
models of a size $k$, i.e. $m(k)=0$. For $m(k)>0$,
we assume all $m(k)$ admissible models of a given size $k$ to be
equally likely, that is, conditionally on the model size $|M|=k$,
$$
P(M \bigl|\; |M|=k)=m(k)^{-1},
$$
where recall that $m(r)=1$ and hence $P(M \bigl|\; |M|=r)=1$.
To complete the prior, for any given model $M$
we assume the normal prior on its unknown coefficients $\bbeta_M$:
$\bbeta_M=\bbeta|M \sim N_p({\bf 0}, \gamma \sigma^2 (D_M X'XD_M)^{+})$ which is
the well-known conventional $g$-prior of Zellner (1986).

For the proposed hierarchical prior, straightforward Bayesian calculus yields
the posterior probability of a model $M$:
\begin{equation}
P(M|\by) \propto\
\pi(|M|) m(|M|)^{-1}
(1+\gamma)^{-\frac{|M|}{2}}
\times \exp\left\{\frac{\gamma}{\gamma+1}\frac{{\bf y}'XD_M(D_M X'XD_M)^{+}D_MX'\by}{2\sigma^2}\right\}, \label{eq:post}
\end{equation}
where we set $\pi(|M|)m(|M|)^{-1}=0$ if $m(|M|)=0$.
Finding the most likely model leads therefore to the following maximum {\em
a posteriori} (MAP) model selection criterion:
$$
\max_M \left\{\by'XD_M(D_M X'XD_M)^{+}D_MX'\by
+2\sigma^2(1+1/\gamma)\ln\left\{m(|M|)^{-1}\pi(|M|)(1+\gamma)^{-\frac{|M|}{2}}\right\} \right\}
$$
or, equivalently,
\begin{equation}
\min_M \left\{||\by-X\hat{\bbeta}_M||^2+2\sigma^2(1+1/\gamma)\ln\left\{m(|M|)\pi(|M|)^{-1}(1+\gamma)^{\frac{|M|}{2}}\right\} \right\} \label{eq:MAP1}
\end{equation}
which is of the general type (\ref{eq:pmle}) with the complexity penalty
\begin{equation} \label{eq:pen}
Pen(|M|)=2\sigma^2(1+1/\gamma)\ln\left\{m(|M|)\pi(|M|)^{-1}(1+\gamma)^{\frac{|M|}{2}}\right\}
\end{equation}
A specific form of the penalty (\ref{eq:pen}) is defined by the choice of
the prior $\pi(\cdot)$ on the model size.

\section{Main results} \label{sec:main}
\subsection{Risk bounds} \label{subsec:bounds}
Denote the set of all admissible models by $\Omega$.
For a given model size $1 \leq p_0 \leq r$, define the set
of all admissible models ${\cal M}_{p_0}$ of size $p_0$, that is,
${\cal M}_{p_0}=\{M \in \Omega: |M|=p_0\}$ and $card({\cal M}_{p_0})=m(p_0)$.
Obviously, if a model $M \in {\cal M}_{p_0}$, the corresponding coefficients vector
$\bbeta_M \in \mathbb{R}^p$ has $p_0$ nonzero components, where
$\beta_{Mj} \neq 0$ iff $X_j$ is included in $M$.
Let $\Mp=\bigcup_{k=0}^{p_0} {\cal M}_k$ be the set of all
admissible models with at most $p_0$ predictors.
Following our arguments from Section \ref{sec:map}, $\Omega$ can be essentially
reduced to $\mathfrak{M}_r$.

In this section we derive the upper and lower bounds for the maximal risk of
the proposed MAP model selector (\ref{eq:MAP1}) over $\Mp$.

\begin{theorem}[upper bound] \label{th:upper}
Let $\hat{M}$ be the solution of (\ref{eq:MAP1}) and $\hat{\bbeta}_{\hat M}$
be the corresponding least squares estimator of its coefficients. Define
$c(\gamma)=8(\gamma+3/4)^2>9/2$ and assume that for some constant $c>0$,
\begin{equation} \label{eq:assump}
\min\{m(k)^{-c},e^{-ck}\} \leq \pi(k) \leq m(k)e^{-c(\gamma)k}
\end{equation}
for all $k=1,\ldots,r$ such that $m(k)>0$.

Then, there exists a constant $C(\gamma)>0$ depending only on $\gamma$ such that
\begin{equation}
\sup_{\bbeta_M: M \in \Mp}E||X\hat{\bbeta}_{\hat M}-X\bbeta_M||^2
\leq C(\gamma) \sigma^2 \min\left\{\max\left(p_0,\ln m(p_0)\right),r\right\}
\label{eq:upper}
\end{equation}
simultaneously for all $1 \leq p_0 \leq r$.
\end{theorem}
Under the conditions (\ref{eq:assump}) on the prior $\pi(k)$ in
Theorem \ref{th:upper}, the corresponding penalty (\ref{eq:pen}) satisfies
\begin{eqnarray} C_1(\gamma) \sigma^2 k & \leq Pen(k) & \leq
C_2(\gamma)\sigma^2 \max\left\{(c+1)\ln m(k)+
\frac{k}{2}\ln(1+\gamma)~;~\ln m(k) +k(c+0.5\ln(1+\gamma))\right\}
 \nonumber \\
&  & \leq C_3(\gamma)\sigma^2\max\left\{\ln m(k) ,k\right\},\;\;\;k=1,\ldots,r \label{eq:penbounds}
\end{eqnarray}
for some positive constants $C_1(\gamma), C_2(\gamma)$ and
$C_3(\gamma)$.

To assess the accuracy of the established upper bound for the quadratic risk of
the proposed MAP estimator, we derive the
lower bound for the minimax risk of estimating $X\bbeta$ in (\ref{eq:model}).

The $l_0$ quasi-norm $||\bbeta||_0$ of a vector $\bbeta$ is defined as the number of its nonzero entries.
For any given $k=1,\ldots,r$, let $\phi_{min}[k]$ and $\phi_{max}[k]$
be the $k$-sparse minimal and maximal eigenvalues of the design defined as
$$
\phi_{min}[k]=\min_{\bbeta: 1 \leq ||\bbeta||_0 \leq k}
\frac{||X\bbeta||^2}{||\bbeta||^2},
$$
$$
\phi_{max}[k]=\max_{\bbeta: 1 \leq ||\bbeta||_0 \leq k}
\frac{||X\bbeta||^2}{||\bbeta||^2}
$$
In fact, $\phi_{min}[k]$ and $\phi_{max}[k]$ are respectively
the minimal and maximal eigenvalues of all $k \times k$ submatrices
of the matrix $X'X$ generated by any $k$ columns of $X$.
Let $\tau[k]=\phi_{min}[k]/\phi_{max}[k],\;k=1,\ldots,r$.
By the definition, $\tau[k]$ is a non-increasing function of $k$.
Obviously, $\tau[k] \leq 1$ and for the orthogonal design
the equality holds for all $k$.

\begin{theorem}[minimax lower bound] \label{th:lower}
Consider the model (\ref{eq:model}) and let $1 \leq p_0 \leq r$.
There exists a universal constant $C>0$ such that
\begin{equation}
\label{eq:lower}
\inf_{\tilde{\bf y}} \sup_{\bbeta_M: M \in \Mp}
\mathbb{E}||\tilde{\bf y}-X\bbeta_M||^2 \geq
\left\{
\begin{array}{ll}
C \sigma^2
\max\left\{\tau[2p_0]\frac{\ln m(p_0) }{\max(1,\ln p_0)},\tau[p_0]p_0\right\},&
\; 1 \leq p_0 \leq r/2  \\
C \sigma^2 \tau[p_0]\; r,  &\; r/2 \leq p_0 \leq r
\end{array}
\right.
\end{equation}
where the infimum is taken over all estimates $\tilde{\bf y}$ of the mean vector
$X\bbeta$.
\end{theorem}

In some particular cases, e.g. for complete variable selection,
the general minimax lower bounds established in Theorem \ref{th:lower} can
be improved  by removing
the $\max(1,\ln p_0)$-term in (\ref{eq:lower}) (Abramovich \& Grinshtein,
2010; Rigollet \& Tsybakov, 2011; Raskutti, Wainwright \& Yu, 2012).
Whether this additional log-term can be removed in the general
case remains so far a conjecture.
Note however that if $\ln m(p_0)=O(p_0)$
(in particular, for the ordered variable selection, where $m(p_0)=1$),
the dominating term in both bounds (\ref{eq:upper}) and (\ref{eq:lower}) is
anyway $p_0$.

The upper bound (\ref{eq:upper}) holds for any design matrix $X$ of
rank $r$, while the minimax lower bound (\ref{eq:lower}) depends on
$X$ but only through the sparse eigenvalues ratios. Finally note
that the structural constraints are manifested in the upper and
lower bounds only through $m(p_0)$ -- the number of admissible
models of size $p_0$.

\subsection{Asymptotic adaptive minimaxity of the MAP estimator} \label{subsec:minmax}
The established upper and lower risk bounds (\ref{eq:upper}),
(\ref{eq:lower}) in the previous Section \ref{subsec:bounds} is the
key for investigating the asymptotic minimaxity of the proposed MAP
estimator, where the number of possible predictors $p=p_n$ may
increase with the sample size $n$.
One can view such a setup as a series of
projections of the vector $X\bbeta$ on the expanding span of
predictors. In particular, it may be $p_n > n$ or even $p_n \gg n$.
Thus, formally, we consider now a {\em sequence} of design matrices
$X_{p,n}$, where $r_n=rank(X_{n,p}) \rightarrow \infty$. For simplicity of exposition, hereafter, we omit
the index $n$. Similarly, there are sequences of the coefficient
vectors $\bbeta_p$ and priors $\pi_p(\cdot)$. In these notations the
original model (\ref{eq:model}) is transformed into a sequence of
models
\begin{equation}
\by=X_p \bbeta_p + \bepsilon, \label{eq:model1}
\end{equation}
where $rank(X_p)=r$ and any $r$ columns of $X_p$ are linearly
independent (hence, $\tau_p[r]>0$),
$\bepsilon \sim N({\bf 0},\sigma^2 I_n)$ and
the noise variance $\sigma^2$ does not depend on $n$ and $p$.

We consider the {\em nearly-orthogonal} design, where the sequence
of sparse eigenvalues ratios $\tau_p[r]$ is bounded away from zero.
Nearly-orthogonality means that there are no ``too strong'' linear
relationships within any set of $r$ columns of the design matrix
$X_p$. Evidently, in this case $p$ cannot be ``too large'' relative
to $r$ and, therefore, to $n$. Indeed, Abramovich \& Grinshtein
(2010) showed that nearly-orthogonality of a design necessarily
implies $p=O(r)$ and, thus, $p=O(n)$. In this case,
$$
\max\left(\ln m(p_0),p_0\right)~ \leq ~\max\left(\ln {p \choose p_0},p_0\right)~ \leq~ p_0(\ln(p/p_0)+1)
~\leq~ p=O(r)
$$
and the following two corollaries are immediate consequences of Theorems
\ref{th:upper} and \ref{th:lower}:
\begin{corollary}[bounds for the minimax risk] \label{cor:minimax}
Let the design be nearly-orthogonal.
There exist two constants $0<C_1 \leq C_2<\infty$ such that for all
sufficiently large $r$,
$$
C_1 \sigma^2 \max\left\{\frac{\ln m(p_0)}{\max(1,\ln p_0)},p_0\right\}
\leq \inf_{\tilde{\by}}\sup_{\bbeta_M: M \in \Mp}E||\tilde{\by}-X_p\bbeta_M||^2
\leq C_2 \sigma^2 \max\left\{\ln m(p_0),p_0\right\}
$$
for all $1 \leq p_0 \leq r$.
\end{corollary}

\begin{corollary}[asymptotic adaptive minimaxity of the MAP estimator]
\label{cor:map} Consider the nearly-orthogonal design and assume
that for $m(k)>0$ the prior $\pi(k)$ satisfies
$$
\min\{m(k)^{-c},e^{-ck}\} \leq \pi(k) \leq m(k)e^{-c(\gamma)k},\;\;\;k=1,\ldots,r
$$
for some $c>0$ and $c(\gamma)$ defined in Theorem \ref{th:upper}. Then the corresponding MAP estimator
(\ref{eq:MAP1}) is asymptotically at least nearly-minimax (up to a $\ln p_0$-factor)
simultaneously over all $\Mp,\;1 \leq p_0 \leq r$.
\end{corollary}

The above general results depend on the asymptotic behavior of
$m(p_0)$ as a function of $p_0$. Similar to Birg\'e \& Massart (2007) we
consider the following three typical cases:

\medskip
\noindent {\bf 1. Small numbers of admissible models}:
$\ln m(k)=O(k),\;k=1,\ldots,r-1$ (recall however that $m(k) \leq {p \choose k}$).
\newline
In particular, for ordered variable selection, $m(k)=1$ though ``small''
numbers allow also even exponential growth of $m(k)$ for small and moderate
$k$.

For this case, Corollary \ref{cor:minimax} and Theorem \ref{th:upper} imply:
\begin{corollary} \label{cor:small}
Consider the nearly-orthogonal design and let $\ln m(k)=O(k),\;k=1,\ldots,r-1$.
As $r$ increases,
\begin{enumerate}
\item
\be \label{eq:msesmall}
\inf_{\tilde{\by}}\sup_{\bbeta_M: M \in \Mp}  E||\tilde{\by}-X_p\bbeta_M||^2
\asymp \sigma^2 p_0
\ee
for all $p_0=1,\ldots,r$.
\item For $m(k)>0$ assume that
$e^{-ck} \leq \pi(k) \leq m(k)e^{-c(\gamma)k},\;k=1,\ldots,r$
for some $c>0$.
The resulting MAP estimator (\ref{eq:MAP1}) attains then the minimax rates
simultaneously over all $\Mp, 1 \leq p_0 \leq r$.
\end{enumerate}
\end{corollary}
Note that $\sigma^2 p_0$ is the risk of (unbiased) least squares estimation
of $X\bbeta_{M_0}$ for a true model $M_0$ of size $p_0$ in (\ref{eq:model}).
Corollary \ref{cor:small} therefore verifies that when the number of
admissible models is small, there is essentially no extra price for model
selection.

It follows from (\ref{eq:pen}) that priors satisfying the conditions of
Corollary \ref{cor:small} lead to the AIC-type penalties of the form
$Pen(k) \sim 2C(\gamma)\sigma^2 k$ for some $C(\gamma)>1$.

\medskip
\noindent {\bf 2. Complete variable selection}:
$m(k)={p \choose k},\;k=1,\ldots,r-1$.
\newline
As we have mentioned above,
from the already known results of Abramovich \& Grinshtein (2010),
Rigollet \& Tsybakov (2011) and Raskutti, Wainwright \& Yu (2012)
the $\max(1,\ln p_0)$-term in the lower minimax risk bound can be removed
in this case:
\begin{corollary} \label{cor:large}
Consider complete variable selection for the nearly-orthogonal design.
As $r$ increases,
\begin{enumerate}
\item
$$
\inf_{\tilde{\by}}\sup_{\bbeta_M: M \in \Mp} E||\tilde{\by}-X_p\bbeta_M||^2
\asymp \sigma^2 p_0 (\ln(p/p_0)+1)
$$
for all $p_0=1,\ldots,r$.
\item Assume that
$\left(\frac{k}{pe}\right)^{ck} \leq \pi(k) \leq \left(\frac{p}{ke^{c(\gamma)}}\right)^k,\;k=1,\ldots,r-1$ and $e^{-cr} \leq \pi(r) \leq e^{-c(\gamma)r}$
for some $c>c(\gamma)$.
The resulting MAP estimator (\ref{eq:MAP1}) attains then the minimax rates
simultaneously over all $\Mp, 1 \leq p_0 \leq r$.
\end{enumerate}
\end{corollary}
The upper bound on $\pi(k)$ in Corollary \ref{cor:large}
trivially holds for all $k \leq pe^{-c(\gamma)}$.

Corollary \ref{cor:large} shows that for complete variable selection,
model selection yields an additional multiplicative factor of $\ln(p/p_0)$
to the risk $\sigma^2p_0$ of estimating $X\bbeta_{M_0}$ for a true model $M_0$
of size $p_0$ in (\ref{eq:model}).

The conditions on the prior of Corollary \ref{cor:large} hold, for example,
for the truncated geometric prior $\pi(k) \propto q^k,\;k=1,\ldots,r$ for some
$0 < q < 1$, and
the corresponding penalties in (\ref{eq:pen}) are of
the $2k\ln(p/k)$-type, where
$Pen(k)=2 C(\gamma) \sigma^2 k(\ln(p/k)+1)(1+o(1))$ for some $C(\gamma)>1$.

\medskip
\noindent {\bf 3. Intermediate case}:
$k=o(\ln m(k))$ and $m(k) <{p \choose k} ,\;k=1,\ldots,r-1$.
\newline
A practically interesting example of the intermediate case is hierarchical model
selection with paired interactions mentioned in Section \ref{sec:intr}:

\medskip
\noindent {\em Example: hierarchical model selection with paired interactions.}
\newline
Consider model selection in regression with $K$ main predictors and their
paired interactions.
The overall number of possible predictors $p$ in
(\ref{eq:model}) is therefore $p=K+{K \choose 2}=K(K+1)/2$.
However, an interaction can be included in the model only together with the
corresponding main effects. Obviously, $m(k) < {p \choose k},\;k=1,\ldots,r-1$.
On the other hand, this is an example with ``large'' numbers of
admissible models, where $\ln m(k) \geq c k \ln(p/k)$ for some
$0 < c < 1$.
Indeed, for models of sizes one and two only
main effects can be selected, and the numbers of admissible models are
$m(1)=K$ and $m(2)={K \choose 2}$ respectively.
One can trivially verify that in both cases
$\ln m(k) \geq c k \ln(p/k),\;k=1,2$ for some
positive constant $c<1$. For $k \geq 3$,
we have the following lemma:
\begin{lemma} \label{lem:lemma1} For all $3 \leq k \leq r-1$,
$$
\ln m(k) \geq \left\lfloor\frac{k}{3}\right\rfloor\ln(p/k)
$$
\end{lemma}

\medskip

For the intermediate case there is the $\ln p_0$-gap between the upper bound for
the quadratic risk of the MAP estimator in Theorem \ref{th:upper}
and the minimax lower bound in Theorem \ref{th:lower}. So far, we
can claim that if $m(k)^{-c} \leq \pi(k) \leq
m(k)e^{-c(\gamma)k)},\;k=1,\ldots,r$ and the design is
nearly-orthogonal, the resulting MAP estimator is asymptotically at
least nearly-minimax (up to $\ln p_0$-factor) over all $\Mp,\; 1
\leq p_0 \leq r$ and can only conjecture that, similar to the
complete variable selection, this log-factor can be removed in this
case as well.

Similar to complete variable selection, when the numbers of
admissible models for the intermediate case are ``large'' (e.g., hierarchical model selection
with main effects and paired interactions considered above),
the conditions on the prior in Corollary \ref{cor:map} are satisfied for
the truncated geometric prior $\pi(k) \propto q^k, k=1,\ldots,r$ for some
$0 < q < 1$ corresponding to complexity penalties of $2k\ln(p/k)$-type.
\medskip
\noindent

Note also that for the nearly-orthogonal design,
$||X_p\hat{\bbeta}_{p\hat M}-X_p\bbeta_p|| \asymp ||\hat{\bbeta}_{p\hat M}-\bbeta_p||$
and all the results of this section for estimating the
mean vector $X_p\bbeta_M$ can be therefore
straightforwardly applied
to estimating the regression coefficients $\bbeta_M$.

\section{Simulation study} \label{sec:simul}
We conducted a short simulation study to demonstrate the performance of
the proposed MAP model selection procedure. We considered polynomial regression which is
an example of the ordered variable selection (see Section \ref{sec:intr}):
$$
y_i=\beta_0+\beta_1 x_i + \ldots + \beta_k x_i^k + \epsilon_i, \quad i=1,\ldots,n,
$$
where $0 \leq k \leq n-1$ is the polynomial degree to be selected, and
$\epsilon_i \sim {\cal N}(0,\sigma^2)$ with the known variance $\sigma^2$
and independent.
In this case obviously $m(k)=1$ for all $k=0,\ldots,n-1$.
An example of a prior satisfying the conditions of Corollary \ref{cor:small} is
the (truncated) geometric prior $Geom(1-q)$, where $\pi(k)=(1-q)q^k/(1-q^n)
\propto q^k,\;
k=0,\ldots,n-1$ for some $0<q<1$. The corresponding complexity penalty
(\ref{eq:pen}) is the AIC-type linear penalty
\begin{equation} \label{eq:pengeom}
Pen(k)=2\sigma^2(1+1/\gamma)\ln\left(q^{-1}\sqrt{1+\gamma}\right)k
\end{equation}

\subsection{Estimation of parameters} \label{subsec:parameter}
To apply the developed MAP model selection procedure we need to
specify the prior parameters $\gamma$ and $q$ in (\ref{eq:pengeom}).
They are rarely known {\em a priori} in practice and usually should
be estimated from the data.

Let $X_{n \times n}$ be the design matrix of the saturated model, that is,
$X_{ij}=x_i^{j-1},\;i=1,\ldots,n; j=0,\ldots,n-1$. For a given $k$, define the
corresponding diagonal indicator matrix $D_k=diag(\bd_k)$, where $d_{kj}=1,\;j=1,\ldots,k$
and zero otherwise (see Section \ref{sec:map}).
Based on the proposed Bayesian model from Section \ref{sec:map},
$$
\by|k \sim {\cal N}({\bf 0},\sigma^2(I+\gamma H_k)),
$$
where $H_k=X D_k(D_k X'XD_k)^+D_kX'$, and straightforward calculus yields
the following marginal likelihood of the
observed data $\by$:
\begin{equation*}
\begin{split}
L(\by;\gamma,q)&=\sum_{k=0}^{n-1} \frac{1}{\sqrt{|I+\gamma H_k|}}
\exp\left\{-\frac{\by'(I+\gamma H_k)^{-1}\by}{2\sigma^2}\right\} \pi(k)\\
& \propto
\sum_{k=0}^{n-1}\frac{1}{(1+\gamma)^{k/2}} \exp\left\{\frac{\gamma}{\gamma+1}
\frac{\by'H_k\by}
{2\sigma^2}\right\} \frac{(1-q)q^k}{1-q^n}
\end{split}
\end{equation*}

The MLEs for $\gamma$ and $q$ can be obtained (numerically) by the EM algorithm.
Regard $k$ as a ``missing'' data and define the corresponding latent indicator
variables $u_j=\delta_{jk},\;j=0,\ldots,n-1$.
The complete log-likelihood for the
``augmented'' data $(\by,\bu)$, up to an additive constant, is then
$$
l(\by,\bu;\gamma,q)=\frac{\gamma}{\gamma+1} \sum_{k=0}^{n-1}u_k
\frac{\by'H_k\by} {2\sigma^2} -
\frac{\ln(1+\gamma)}{2}\sum_{k=0}^{n-1}u_k k + \ln q \sum_{k=0}^{n-1}u_k k
+\ln(1-q)-\ln(1-q^n)
$$

On the E-step at the $h$-th iteration we compute the conditional expectation
\begin{equation}
\begin{split}
\hat{l}^{[h]}&=E\left(l(\by,\bu;\gamma,q)|\by,\gamma^{[h]},q^{[h]}\right)\\
&=\frac{\gamma}{\gamma^+1} \sum_{k=0}^{n-1}\hat{u}^{[h]}_k
\frac{\by'H_k\by} {2\sigma^2} -
\frac{\ln(1+\gamma)}{2}\sum_{k=0}^{n-1}\hat{u}^{[h]}_k k + \ln q \sum_{k=0}^{n-1}\hat{u}^{[h]}_k k
+\ln(1-q)-\ln(1-q^n), \label{eq:comlike}
\end{split}
\end{equation}
where
$$
\hat{u}^{[h]}_k=E(u_k|\by,\gamma^{[h]},q^{[h]})=
\frac{\exp\left\{\frac{\gamma^{[h]}}{\gamma^{[h]}+1} \frac{\by'H_k\by}{2\sigma^2}
\right\} q^{[h]k}}
{\sum_{j=0}^{n-1}\exp\left\{\frac{\gamma^{[h]}}{\gamma^{[h]}+1} \frac{\by'H_j\by}{2\sigma^2}
\right\} q^{[h]j}}
$$

At the M-step we maximize $\hat{l}^{[h]}$ w.r.t. $\gamma$ and $q$ to get
$$
\hat{\gamma}^{[h+1]}=\left(\frac{\sum_{k=0}^{n-1} \hat{u}^{[h]}_k \by'H_k\by}
{\sigma^2\sum_{k=0}^{n-1} \hat{u}^{[h]}_k k}-1\right)_+
$$
There is no closed form solution for $\hat{q}^{[h+1]}$. However,
replacing the truncated geometric distribution by a usual one and
ignoring thus the last term $\ln(1-q^n)$ in the RHS of (\ref{eq:comlike}),
implies a good approximation of $\hat{q}^{[h+1]}$ for large $n$:
$$
\hat{q}^{[h+1]}=\frac{\sum_{k=0}^{n-1} \hat{u}^{[h]}_k k}
{1+\sum_{k=0}^{n-1} \hat{u}^{[h]}_k k}
$$

\subsection{The results} \label{subsec:results}
We used two test functions: a fifth degree polynomial
$g_1(x)=(x+0.1)(x-0.2)(x-0.4)(x-0.8)(x-1.1)$ and the Doppler
function $g_2(x)=\sqrt{x(1-x)}\sin(2\pi\cdot1.05/(x+0.05))$ (see Donoho \&
Johnstone, 1994) which does not have a sparse polynomial approximation.
Both functions were then normalized to have unit $L_2[0,1]$-norms.
The data were generated according to the model
$$
y_i=g_{1,2}(i/n)+\epsilon_i,\quad i=1,\ldots,n
$$
for $n=100$, where $\epsilon_i \sim {\cal N}(0,\sigma^2)$ and independent.
The noise variance $\sigma^2$ was chosen to ensure the signal-to-noise ratio
SNR at levels $3, 5$ and $7$ and was assumed to be known.
The number of replications was 100.

The proposed MAP model selector results in a model selection procedure with a linear type
penalty of the form $Pen(k)=2\sigma^2\lambda k$ with
$\lambda_{MAP}=(1+1/\gamma)\ln(q^{-1}\sqrt{1+\gamma})$ (see (\ref{eq:pengeom})). We compared it
with two other well-known model selection procedures
with linear penalties, namely,
AIC ($\lambda_{AIC}=1$) and RIC ($\lambda_{RIC}=\ln p$) (see Section \ref{sec:intr}). In our case, $\lambda_{RIC}=\ln 100=4.6$.

Table \ref{tab:mse} summarizes mean squared errors averaged over 100
replications (AMSE). We also present the average
polynomial degrees selected by the three methods for approximating the true
response functions.

\begin{table}[htb]
\begin{center}
\caption{AMSEs and polynomial degrees (in parentheses) averaged over
100 replications for three estimators}
\begin{tabular}{c|ccc|ccc}
    & \multicolumn{3}{c|}{$g_1$} &\multicolumn{3}{c}{$g_2$}\\
\hline
SNR &   MAP  &  AIC   &  RIC   &  MAP   &  AIC   &  RIC  \\
\hline
3   & 0.661  & 0.813  & 0.652  & 28.844 & 27.818 & 35.188 \\
    &(5.01)  & (5.50) & (5.00) &(29.07) & (33.66)& (18.11)\\
5   & 0.238  & 0.293  & 0.235  & 24.677 & 23.450 & 27.056 \\
    &(5.01)  & (5.50) & (5.00) &(44.13) & (53.48)& (29.40)\\
7   & 0.121  & 0.149  & 0.120  & 22.995 & 22.698 & 26.036 \\
    &(5.01)  & (5.50) & (5.00) &(52.72) & (58.57)& (32.21)\\
\end{tabular}
\label{tab:mse}
\end{center}
\end{table}

As expected, a more conservative RIC tends to include less terms in the model
and outperforms AIC for a polynomial
$g_1$, while the latter is superior for $g_2$, where a high
order polynomial approximation is required.
The MAP estimator with estimated $\gamma$
and $q$ yields a data-driven
$\lambda$ and is adaptive to the unknown polynomial degree -- it behaves very
similar to RIC when it is low and to AIC when it is high.

\section{Concluding remarks} \label{sec:tam}
In this paper we considered model selection in linear regression under
general structural constraints and extended the existing results for the
complete
variable selection. In particular, we utilized a Bayesian MAP model selection
procedure of Abramovich \& Grinshtein (2010) and modified it correspondingly.
From a frequentist view, the resulting MAP model selector is
a penalized least squares estimator with a complexity penalty associated with
a prior on the model size which is adaptive to the structural constraints.
In fact, the proposed
Bayesian approach can be used as a tool for generating a wide class of
penalized least squares estimators with various complexity penalties.

We established the general upper bound for the quadratic risk of the MAP estimator
over a set of admissible models of a given size and the lower bound for the
corresponding minimax risk.
Based on these results, we showed that for the
nearly-orthogonal design, the MAP estimator is asymptotically at least nearly
minimax (up to a log-factor)
for the entire range of sparse and dense models. Moreover, when the
numbers of admissible models are ``small'' or, on the opposite, for the case of complete
variable selection, it achieves the exact
minimax rates. The corresponding MAP model selection procedures lead respectively to
AIC-type and $2k\ln(p/k)$-type criteria.
Whether these results on asymptotic minimaxity are true for the intermediate
case remains so far a conjecture.

There are also other challenges for future research. The
assumption of nearly-orthogonal design used in investigating
asymptotic minimaxity of MAP estimators in Section
\ref{subsec:minmax} typically does not hold for $p \gg n$ setup due
to the multicollinearity phenomenon. The analysis of multicollinear
design, where the sequence of sparse eigenvalues ratios $\tau_p[r]$ may
tend to zero as $p$ increases is much more delicate. In this case
there is a gap (in addition to a log-factor) between the rates in
the upper and lower bounds (\ref{eq:upper}) and (\ref{eq:lower}).
Unlike model identification or coefficients estimation, where
multicollinearity is a ``curse'', it may essentially become a
``blessing'' for estimating the mean vector allowing one to exploit
correlations between predictors to reduce the size of a model
(hence, to decrease the variance) without paying much extra price in
the bias term. Interestingly, a similar phenomenon also occurs in a
testing setup (e.g., Hall \& Jin, 2010). Abramovich \& Grinshtein
(2010) investigated the complete variable selection for
multicollinear design and showed that under certain additional
assumptions on the design and the regression coefficients, the MAP
estimator corresponding to the  $2k\ln(p/k)$-type complexity penalty
remains asymptotically rate-optimal (in the minimax sense) even for
this case. Whether it is true and what are the additional conditions
in the presence of structural constraints is a challenging topic for
future research.

Computational issues are another important problem. When the numbers
of admissible models are large (e.g., complete variable selection or
hierarchical model selection with interactions), minimizing
(\ref{eq:pmle}) (and (\ref{eq:MAP1}) in particular) requires
generally an NP-hard combinatorial search. During the last decade
there have been substantial efforts to develop various approximated
algorithms for solving (\ref{eq:pmle}) that are computationally
feasible for high-dimensional data (see, e.g. Tropp \& Wright, 2010
for a survey and references therein). The common remedies involve
either greedy algorithms (e.g., forward selection, matching pursuit)
approximating the global solution by a stepwise sequence of local
ones, or convex relaxation methods replacing the original
combinatorial problem by a related convex program (e.g., Lasso and
Dantzig selector for linear penalties). Abramovich \& Grinshtein
(2010) proposed to utilize the developed Bayesian formalism for
solving (\ref{eq:MAP1}) by using a stochastic search variable
selection (SSVS) techniques originated in George \& McCullogh (1993,
1997). The underlying idea of SSVS is based on generating a sequence
of models from the posterior distribution $P(M|{\bf y})$ in
(\ref{eq:post}). The key point is that the relevant models with the
highest posterior probabilities will appear most frequently and can
be identified even for a generated sample of a relatively small size
avoiding computations of the entire posterior distribution. However,
most of the above approaches have been developed and studied for
complete variable selection. Their adaptation to minimization of
(\ref{eq:MAP1}) subject to the additional structural constraints
while remaining computationally feasible should depend on the
specific type of constraints at hand. In particular, for somewhat
different priors, Chipman (1996) and Farcomeni (2010) considered
SVSS for hierarchical model selection in regression with paired
interactions (see the example in Section \ref{sec:main}) and for
model selection with factor predictors, where the corresponding
dummy variables are all included or excluded. Bien, Taylor \&
Tibshirani (2012) modified Lasso for hierarchical paired
interactions (see also references therein). However, to the best of
our knowledge, there are no theoretical results on the optimality of
the resulting estimators.

Finally, we should note that the obtained theoretical results assume
that the noise variance $\sigma^2$ is known which is rarely the case
in practical applications. One can estimate $\sigma^2$ from the data
with the additional tuning parameters of the MAP procedure (see
Section \ref{subsec:parameter}). Alternatively, he can follow the
fully Bayesian approach and impose some prior distribution on it
(see, e.g., Chipman, 1996; Chipman, George \& McCulloch, 2001 and
Farcomeni, 2010).

\medskip
\noindent
{\bf Acknowledgments}. Both authors were supported by the Israel
Science Foundation grant ISF-248/08.

\section{Appendix}
Throughout the proofs we use $C$ to denote a generic positive constant, not
necessarily the same each time it is used, even within a single equation.

\subsection{Proof of Theorem \ref{th:upper}}
We first show that under the conditions on a prior in Theorem \ref{th:upper},
the corresponding penalty $Pen(k)$ in (\ref{eq:pen}) belongs to the
class of penalties considered in Birg\'e \& Massart (2001) and then use
their Theorem 2 to establish the general upper bound for the quadratic risk
of the MAP estimator (\ref{eq:MAP1}).

Define
\be \label{eq:Lk}
L_k=\frac{1}{k}\ln\left(m(k)\pi^{-1}(k)\right),\quad k=1,\ldots,r,
\ee
where under the conditions on the prior $\pi(\cdot)$ in Theorem
\ref{th:upper}, $L_k \geq c(\gamma)$.
In terms of $L_k$ the complexity penalty (\ref{eq:pen}) is
$Pen(k)=\sigma^2(1+1/\gamma)k(2L_k+\ln(1+\gamma))$.

In our notations the conditions (3.3) and (3.4) on $L_k$ in Theorem 2
of Birg\'e \& Massart (2001) correspond respectively to
\be \label{eq:cond1}
\sum_{k=1}^r m(k)e^{-k L_k} < c
\ee
and
\be \label{eq:cond2}
(1+1/\gamma)(2L_k+\ln(1+\gamma)) \geq C(1+\sqrt{2L_k})^2,\quad k=1,\ldots,r
\ee
for some $C>1$.

The condition (\ref{eq:cond1}) follows immediately from the definition of $L_k$:
$$
\sum_{k=1}^r m(k)e^{-k L_k}=\sum_{k=1}^r \pi(k)=1-\pi(0)<1
$$

Consider now (\ref{eq:cond2}) which is equivalent
to the inequality
\be \label{eq:cond21}
2(1+1/\gamma-C)L_k-2C\sqrt{2L_k}+(1+1/\gamma)\ln(1+\gamma)-C \geq 0
\ee
Repeating the calculus in the proof of Theorem 1 of Abramovich, Grinshtein
\& Pensky (2007), one verifies that with the upper bound on the prior $\pi(k)$
in (\ref{eq:assump}), for $C=1+1/(2\gamma)$,
(\ref{eq:cond21}) and therefore (\ref{eq:cond2}) are satisfied for all
$L_k,\;k=1,\ldots,r$.

Given (\ref{eq:cond1})-(\ref{eq:cond2}), Theorem 2 of Birg\'e \& Massart (2001)
yields the following upper bound for the quadratic risk of the MAP estimator
(\ref{eq:MAP1}) of the mean vector $X\bbeta$ in (\ref{eq:model}):
\be
E||X\hat{\bbeta}_{\hat{M}}-X\bbeta||^2  \leq  c_0(\gamma)
\inf_{M \in \mathfrak{M}_r}\left\{||X\bbeta_M-X\bbeta||^2+Pen(|M|)\right\}
+ c_1(\gamma)\sigma^2 \label{eq:upperrisk}
\ee
for some $c_0(\gamma)$ and $c_1(\gamma)$ depending only on $\gamma$, where
$Pen(|M|)$ is given in (\ref{eq:pen}).

Recall that $m(r)=1$. Then, using the lower bound on $\pi(r)$,
(\ref{eq:upperrisk}) and (\ref{eq:penbounds}) imply
\begin{eqnarray}
\sup_{\bbeta_M: M \in \Mp} E||X\hat{\bbeta}_{\hat{M}}-X\bbeta_M||^2
& \leq & \sup_{\bbeta_M: M \in \mathfrak{M}_r} E||X\hat{\bbeta}_{\hat{M}}-X\bbeta_M||^2
~\leq~  c_0(\gamma)Pen(r)+c_1(\gamma)\sigma^2 \nonumber \\
& \leq  & C(\gamma)\sigma^2 r \label{eq:riskr}
\end{eqnarray}
for  all $1 \leq p_0 \leq r$.

On the other hand, for $p_0<r$ from (\ref{eq:upperrisk}) and (\ref{eq:penbounds})
we have
\be
\sup_{\bbeta_M: M \in \Mp}E||X\hat{\bbeta}_{\hat{M}}-X\bbeta_M||^2
~\leq~  c_0(\gamma)Pen(p_0)+c_1(\gamma)\sigma^2~ \leq~  C(\gamma)\sigma^2 \max\left\{\ln m(p_0) ,p_0\right\}
\label{eq:riskp0}
\ee
Combining (\ref{eq:riskr}) and (\ref{eq:riskp0}) completes the proof.
\newline $\Box$

\subsection{Proof of Theorem \ref{th:lower}}
We first show that
\be \label{eq:lbound1}
\inf_{\tilde{\by}} \sup_{\bbeta_M: M \in \Mp}
\mathbb{E}||\tilde{\by}-X\bbeta_M||^2 \geq \tau[p_0]\sigma^2 p_0
\ee
for all $1 \leq p_0 \leq r$.

Consider the original regression model (\ref{eq:model}) and assume that
the true coefficients vector $\bbeta_M \in \Mp$. Define $X_M=XD_M$, where
the $p \times p$ diagonal indicator matrix $D_M$ was defined in Section
\ref{sec:map}.
In the coefficients domain one then has
\be \label{eq:bw}
\bw=\bbeta_M+\bepsilon',
\ee
where $\bw=(X_M'X_M)^{+}X_M'\by$ and $\bepsilon' \sim {\cal N}({\bf 0},\sigma^2(X_M'X_M)^{+})$.

Let $R(\Mp,\sigma^2 (X_M'X_M)^{+})$ be the minimax risk of estimating
$\bbeta_M$ in (\ref{eq:bw}) over $\Mp$, that is,
$$
R(\Mp,\sigma^2 (X_M'X_M)^{+})=\inf_{\tilde{\bbeta}} \sup_{\bbeta_M: M \in \Mp}
\mathbb{E}||\tilde{\bbeta}-\bbeta_M||^2
$$
Evidently,
\be \label{eq:r1}
\inf_{\tilde{\bf y}} \sup_{\bbeta_M: M \in \Mp}
\mathbb{E}||\tilde{\bf y}-X\bbeta_M||^2 \geq
\phi_{min}[p_0]~ R(\Mp,\sigma^2 (X_M'X_M)^{+})
\ee
Consider also the model (\ref{eq:bw}) but with the uncorrelated noise
$\bepsilon'' \sim {\cal N}({\bf 0},\sigma^2\phi^{-1}_{max}[p_0] D_M)$:
\be \label{eq:bw1}
\bw=\bbeta_M+\bepsilon''
\ee
and the corresponding minimax risk $R(\Mp,\sigma^2\phi^{-1}_{max}[p_0] D_M)$.

Since $(X_M'X_M)^{+} \geq \phi^{-1}_{max}[p_0]D_M$ in the usual sense that
$(X_M'X_M)^{+}-\phi^{-1}_{max}[p_0]D_M$ is positive semi-definite,
\be \label{eq:r2}
R(\Mp,\sigma^2 (X_M'X_M)^{+}) \geq R(\Mp,\sigma^2\phi^{-1}_{max}[p_0] D_M)
\ee
(see, e.g., Lemma 4.27 of Johnstone, 2011).

No estimator of $\bbeta_M$ in (\ref{eq:bw1}) obviously cannot outperform the
oracle estimator that knows the true $\bbeta_M \in \Mp$ whose ideal minimal quadratic
risk is $\sum_{j=1}^p \min(\beta^2_{M,j},\sigma^2\phi^{-1}_{max}[p_0])$
(e.g., Donoho \& Johnstone, 1994).
Hence,
\be \label{eq:r3}
R(\Mp,\sigma^2\phi^{-1}_{max}[p_0] D_M) \geq \sup_{\bbeta_M: M \in \Mp}
\sum_{j=1}^p\min(\beta^2_{M,j},\sigma^2 \phi^{-1}_{max}[p_0]) \geq
\sigma^2 \phi^{-1}_{max}[p_0]~p_0
\ee
Combining (\ref{eq:r1}), (\ref{eq:r2}) and (\ref{eq:r3}) implies (\ref{eq:lbound1}).

Consider now  $1 \leq p_0 \leq r/2$ and show that in this case, in addition
to the lower bound (\ref{eq:lbound1}),
$$
\inf_{\tilde{\by}} \sup_{\bbeta_M: M \in \Mp} \mathbb{E}||\tilde{\by}-X\bbeta_M||^2 \geq
C \sigma^2\tau[2p_0]\frac{\ln m(p_0) }{\max(1,\ln p_0)}
$$
Let
\begin{equation}
s^2(p_0)=C \sigma^2\tau[2p_0]\frac{\ln m(p_0) }{\max(1,\ln p_0)} \label{eq:sp0}
\end{equation}
The core of the proof is to find a subset ${\cal B}_{p_0}$ of vectors
$\bbeta_M, M \in \Mp$
and the corresponding subset of mean vectors
${\cal G}_{p_0}=\{{\bf g} \in \mathbb{R}^n: {\bf g}=X\bbeta_M,\;\bbeta_M \in
{\cal B}_{p_0}\}$ such that for any ${\bf g}_1,\;{\bf g}_2 \in {\cal G}_{p_0}$,
$||{\bf g}_1-{\bf g}_2||^2 \geq 4s^2(p_0)$ and the Kullback-Leibler divergence
$K(\mathbb{P}_{\bg_1},\mathbb{P}_{\bg_2})=\frac{||\bg_1-\bg_2||^2}{2\sigma^2} \leq (1/16)
\ln {\rm card}({\cal G}_{p_0})$.
Lemma A.1 Bunea, Tsybakov \& Wegkamp (2007) will imply then that
$s^2(p_0)$ is the minimax lower bound over $\Mp$.

The standard techniques for constructing such sets of vectors
for the complete variable selection setup is based on generalizations of Varshamov-Gilbert
bound (e.g, Abramovich \& Grinshtein, 2010; Rigollet \& Tsybakov, 2011; Raskutti,
Wainwright \& Yu, 2012). Unfortunately, it
cannot be applied when there are additional structural constraints on the set of
admissible models. We utilize instead the recent combinatorial results of
Gutin \& Jones (2012).

Define the subset $\tilde{\cal B}_{p_0}$ of all vectors $\bbeta_M,
M \in {\cal M}_{p_0} \subset \Mp$
that have $p_0$ entries equal to $C_{p_0}$ defined later,
while the remaining  entries are zeros:
$\tilde{\cal B}_{p_0}=\{\bbeta_M: \bbeta_M \in \{\{0,C_{p_0}\}^{p}\},\;
M \in {\cal M}_{p_0}\}$.
Let $\rho(\bbeta_{1M},\bbeta_{2M})=\sum_{j=1}^p \mathbb{I}\{\beta_{1M,j} \neq \beta_{2M,j}\}$
be the Hamming distance between $\bbeta_{1M},\;\bbeta_{2M} \in \tilde{\cal B}_{p_0}$
and define
$$
\rho_{max}=\max_{\bbeta_{1M},\;\bbeta_{2M} \in \tilde{\cal B}_{p_0}}\rho(\bbeta_{1M},\bbeta_{2M}) \;\;\;{\rm and} \;\;\;
\rho_{min}=\min_{(\bbeta_{1M} \ne \bbeta_{2M}) \in \tilde{\cal B}_{p_0}}\rho(\bbeta_{1M},\bbeta_{2M})
$$
Theorem 2 of Gutin \& Jones (2012) ensures that for any constant
$C>2$ there exists a subset ${\cal B}_{p_0} \subset \tilde{\cal
B}_{p_0}$ such that
$$
\frac{\rho_{max}}{\rho_{min}} \leq C \;\;\;\;\;{\rm and}\;\;\;\;\;
\ln{\rm card}({\cal B}_{p_0}) \geq \alpha \ln m(p_0),
$$
where
$$
\alpha=\left\lceil \frac{\ln(p_0/2)}{\ln(C/2)}\right\rceil^{-1}
$$

Consider the corresponding subset of mean vectors ${\cal G}_{p_0}$,
where ${\rm card}({\cal G}_{p_0})={\rm card}({\cal B}_{p_0})$.
For any $\bg_1,\; \bg_2 \in {\cal G}_{p_0}$ and the
associated with them $\bbeta_{1M},\;\bbeta_{2M} \in {\cal B}_{p_0}$ we then have
\begin{equation}
||\bg_1-\bg_2||^2 = ||X(\bbeta_{1M}-\bbeta_{2M})||^2 \geq \phi_{min}[2p_0] \;
||\bbeta_{1M}-\bbeta_{2M}||^2 \geq \phi_{min}[2p_0]C^2_{p_0}\rho_{min}
\label{eq:t4.1}
\end{equation}

On the other hand, by similar arguments, the Kullback-Leibler divergence
satisfies
\begin{equation}
K(\mathbb{P}_{\bg_1},\mathbb{P}_{\bg_2})
\leq
\frac{\phi_{max}[2p_0] C^2_{p_0} \rho(\bbeta_{1M},\bbeta_{2M})}{2\sigma^2}
\leq \frac{\phi_{max}[2p_0] C^2_{p_0} \rho_{max}}{2\sigma^2}
\label{eq:t4.2}
\end{equation}
Set now
$$
C^2_{p_0}=\frac{\sigma^2 \alpha \ln m(p_0)}{8\rho_{\max}\phi_{max}[2p_0]}
$$
Then, (\ref{eq:t4.1}) and (\ref{eq:t4.2}) yield
$||\bg_1-\bg_2||^2 \geq \tau[2p_0]\sigma^2\alpha \ln m(p_0)/(8C)$,
$K(\mathbb{P}_{\bg_1},\mathbb{P}_{\bg_2}) \leq
(1/16)\ln {\rm card}({\cal G}_{p_0})$, and
Lemma A.1 of Bunea, Tsybakov \& Wegkamp (2007) with $s^2(p_0)$ from (\ref{eq:sp0})
completes the proof.
\newline $\Box$

\subsection{Proof of Lemma \ref{lem:lemma1}}
\noindent 1. $3 \leq k \leq \frac{3}{2}K$
\newline
Consider first all models that include any $\lfloor \frac{k}{3}\rfloor$
paired interactions and the corresponding main effects.
Evidently, the size of any such a model is at most $k$. If it is less than $k$
(it happens when the same main effect appears in several interactions),
we complete it to $k$ by adding other main effects from the remaining ones
(one can easily verify that there are enough remaining main effects since
$k-\lfloor \frac{k}{3}\rfloor \leq K$). We then have
$$
m(k) \geq {\frac{K(K-1)}{2} \choose \lfloor\frac{k}{3}\rfloor}
$$
and, therefore, by straightforward calculus,
$$
\ln m(k) \geq \left\lfloor\frac{k}{3}\right\rfloor \ln\left(\frac{K(K-1)}{2\lfloor\frac{k}{3}\rfloor}\right)
\geq \left\lfloor\frac{k}{3}\right\rfloor \ln\left(\frac{K(K+1)}{2k}\right)=
\left\lfloor\frac{k}{3}\right\rfloor \ln(p/k)
$$

\medskip
\noindent 2. $\frac{3}{2}K < k \leq r-1$
\newline
For this case we consider models of size $k$ that include all $K$ main effects
and any $k-K$ paired interactions. Thus,
$$
m(k) \geq {\frac{K(K-1)}{2} \choose k-K}
$$
and
$$
\ln m(k) \geq (k-K)\ln \left(\frac{K(K-1)}{2(k-K)}\right) \geq
\frac{k}{3}\ln\left(\frac{K(K-1)}{2(k-K)}\right)
$$
To complete the proof one can easily verify that
$\frac{K(K-1)}{2(k-K)} \geq \frac{K(K+1)}{2k}$.
\newline $\Box$

\end{document}